\documentclass[letterpaper,12pt]{amsart}
\usepackage[english]{babel}
\usepackage{amsthm}
\usepackage{amsmath}
\usepackage{amssymb}
\usepackage[latin1]{inputenc}
\def\R{\text{$\mathbb{R}$}}

\def\pr{\parallel}

\def\lra{\longrightarrow}

\def\lra{\longrightarrow}

\def\d1#1#2{\frac{d#1}{d#2}}

\def\p1#1#2{\frac{\partial #1}{\partial #2}}
\def\to{t_o}
\def\pr{\parallel}

\def\part{a=\to \leq t_1 \leq \ldots \leq t_{n-1} \leq t_{n}=b}

\def\P{\text{$\mathbb{P}$}}
\def\C{\text{$\mathbb{C}$}}
\def\H{\text{$\mathbb{H}$}}
\newcommand\SO{{\rm SO}}

\newcommand\SU{{\rm SU}}

\newcommand\U{{\rm U}}
\newcommand\SUr[2]{{{\rm S}({\rm U}({#1}) \times {\rm U}({#2}))}}

\newcommand\mf[1]{\mathfrak{#1}}

\def\To{\rm T}
\def\o{\omega}
\def\omr{\omega_{\lambda}}

\def\mr{ M_{\lambda}}
\def\mk{M^{\lambda}}

\def\gd{\mathfrak{g}^*}
\newcommand\sr[1]{M_{#1}}
\title[Reduction and cut spaces]{Some results on Multiplcity-free spaces}
\author[L. Biliotti]{Leonardo Biliotti}
\address{Dipartimento di Matematica, Universit\`a Politecnica delle
Marche, Via Brecce Bianche, 60131, Ancona  Italy}
 \email{biliotti@dipmat.univpm.it}
\begin{document}
\newtheorem{thm}{Theorem}[section]
\newtheorem{prop}[thm]{Proposition}
\newtheorem{lemma}[thm]{Lemma}
\newtheorem{cor}[thm]{Corollary}
\theoremstyle{definition}
\newtheorem{defini}[thm]{Definition}
\newtheorem{notation}[thm]{Notation}
\newtheorem{exe}[thm]{Example}
\newtheorem{conj}[thm]{Conjecture}
\newtheorem{prob}[thm]{Problem}
\theoremstyle{remark}
\newtheorem{rem}[thm]{Remark} 
\begin{abstract}
Let $(M,\omega)$ be a connected symplectic manifold on which a
connected Lie group $G$ acts properly and in a Hamiltonian fashion
with moment map $\mu:M \lra \mf g^*$. Our purpose is investigate 
multiplicity-free actions, giving criteria to decide
a multiplicity  freenes of the action. As an application we give the
complete classification of multiplicity-free actions of 
compact Lie groups acting isometrically and in a Hamiltonian fashion
on Hermitian symmetric spaces of noncompact type. 
Successively we make a connection
between multiplicity-free actions on $M$ and multiplicity-free actions
on the symplectic reduction and
on the symplectic cut, which allow us to give new examples of 
multiplicity-free actions. 
\end{abstract}
\maketitle
\section{Introduction}
Let $(M, \omega)$ be a connected symplectic manifold with a proper
and Hamiltonian action of a connected Lie group $G$ and let $\mu:
M \lra \gd$ be a corresponding moment map. In 1984 Guillemin and
Sternberg \cite{GS}, motivated by geometric quantization,
introduced the notion of multiplicity-free space when the ring
of the $G-$invariant functions on $M$ is commutative with respect
to the Poisson-bracket. The manifold $M$ is called $G$
multiplicity-free space and the $G-$action is called
multiplicity-free.
The term multiplicity-free comes from the representation theory of
Lie groups. 

A unitary representation of a Lie group $G$ on a
Hilbert space $\H$ is said to be multiplicity-free if each
irreducible representation of $G$ occurs with multiplicity zero or one in $\H$.
The relationships between the two definitions comes via the theory of geometric
quantization. The condition that a unitary representation of $G$ on $\H$ be
multiplicity-free is equivalent to the condition that the ring of bounded
operators on $\H$ be commutative.

In this paper we investigate multiplicity-free actions, which we also may 
call coisotropic actions, on a symplectic manifold $M$, imposing that $G$ acts
properly and in a Hamiltonian fashion on $M$ and a technical condition,
which is needed for applying the symplectic slice, see \cite{lb} and 
\cite{sl}, and the symplectic stratification of the reduced space given
in \cite{lb}, \cite{sl}, which is the following.

For every $\alpha \in \mf g^*$, $\mf g^*$ is the dual of the Lie
algebra of $G$, $G\alpha$ is a locally closed coadjoint orbit of
$G$. Observe that the condition of a coadjoint orbit being locally
closed is automatic for reductive groups and for their semidirect
products with
vector spaces. There exists an example of a solvable group due to
Mautner \cite{mau} p.$512$, with non-locally closed coadjoint
orbits. 

One of our purpose is to give Equivalence theorem for multiplicity-free action,
which shall allow us to prove that the complete classification
of compact Lie groups acting multiplicity-free
on irreducible Hermitian symmetric spaces of noncompact type follows from one
given in a compact case. 

We also prove a reduction principle for multiplicity-free actions and 
we make a connection between multiplicity-free actions on $M$ and 
multiplicity-free actions on 
the symplectic reduction and on the symplectic cut, mainly in a K\"ahler 
geometry. As an application we give new examples of multiplicity-free actions
on compact K\"ahler manifolds which are not Hermitian symmetric spaces.   

\section{Hamiltonian viewpoint}
Let $M$ be a connected differential manifold equipped with a
non-degenerate closed $2-$form $\omega$. The pair $(M,\omega)$ is called
symplectic manifold. Here we consider a finite-dimensional
connected Lie group acting smoothly and properly on $M$ so that
$g^* \omega=\omega$ for all $g \in G,$ i.e. $G$ acts as a group of
canonical or symplectic diffeomorphism. For $f,g \in
C^{\infty}(M),$ define $\{ f,g \}= \omega (X_f , X_g),$ where
$X_f$ and $X_g$ are the vector fields which is uniquely defined by
$df=i_{X_f} \omega$ and $dg= i_{X_g} \omega.$ It follows that
$(C^{\infty}(M), \{,\})$ is a Poisson algebra.

The $G$-action is called {\em Hamiltonian}, and we said that $G$
acts in a Hamiltonian fashion, if there exists a map
$$
\mu: M \lra \mf{g}^*,
$$
which is called moment map, satisfying:
\begin{enumerate}
\item for each $X\in \mf{g}$ let
\begin{itemize}
\item $\mu^{X}: M \lra \R,\ \mu^{X} (p)= \mu(p) (X),$ the
component of $\mu$ along $X,$ \item $X^\#$ be the vector field on
$M$ generated by the one para\-me\-ter subgroup $\{ \exp (tX):
t \in \R \} \subseteq G$.
\end{itemize}
Then
$$
{\rm d} \mu^{X}= {\rm i}_{X^\#} \omega
$$
i.e. $\mu^{X}$ is a Hamiltonian function for the vector field
$X^\#.$ \item $\mu$ is $G-$equivariant, i.e. $\mu (gp)=Ad^* (g)
(\mu(p)),$ where $Ad^*$ is the coadjoint representation on $\mf g^*$.
\end{enumerate}

Let $x \in M$ and ${\rm d} \mu_x: T_x M \lra T_{\mu(x)} \mf g^*$
be the differential of $\mu$ at $x$. Then
\begin{equation} \label{nucleo}
Ker {\rm d} \mu_x = (T_x G(x))^{\perp_{\omega}}:= \{ v \in T_x M:
\omega (v,w)=0,\ \forall w \in  T_x G(x)\}.
\end{equation}
If we restrict $\mu$ to a $G-$orbit $Gx$, then we have the
homogeneous fibration $\mu: Gx \lra Ad^*(G)\mu(x)$ and the
restriction of the ambient symplectic form $\omega$ on the orbit
$Gx$ is the pullback by the moment map $\mu$ of the symplectic
form on the coadjoint orbit through $\mu(x)$
\begin{equation} \label{KKS}
\omega_{|_{Gx}}= \mu^* (\omega_{ Ad^*(G) \mu(x)})_{|_{Gx}}
\end{equation}
see \cite{lb} p. 211, where $\omega_{G \mu(x)}$ is the
Kirillov-Konstant-Souriau (KKS) symplectic form on the coadjoint
orbit of $\mu(x)$ in $\mf g^*$. This implies the following 
result.
\begin{prop} \label{c1}
A $G$-orbit $Gx$ is a symplectic submanifold of $M$ if and only if
the moment map restricted to $Gx$ into $G\mu(x)$, $\mu_{|Gx}:Gx \longrightarrow G\mu(x)$,
is a covering map.
In particular if $G$ is a compact Lie group then $G_x=G_{\mu(x)}$ and
$\mu_{|_{Gx}}:Gx \longrightarrow G\mu(x)$ is a diffeomorphism onto.
\end{prop}
\begin{proof}
The first affirmation follows immediately from (\ref{KKS}). If $G$
is compact, coadjoint orbits are of the form $G/C(T)$, where
$C(T)$ is the centralizer of the torus $T$. In particular such
orbits are simply connected, form which one may deduce
$G_x=G_{\mu(x)}$.
\end{proof}
\section{Multiplicity-free spaces}
Let $(M,\omega)$ be a connected symplectic manifold and let $G$ be
a connected Lie group acting properly and as a group of symplectic
diffeomorphism on $M$. 
\begin{defini}
The $G$-action is called \emph{multiplicity-free},
$M$ is called a $G$ multiplicity-free space, if the space of
invariant function $C^{\infty}(M)^G$ is a commutative Lie algebra with 
respect the Poisson bracket.
\end{defini}
By the definition follows that if  $K \subset G$ and $M$ is a $K$
multiplicity-free space then $M$ is a $G$ multiplicity-free space
as well.

The multiplicity-free actions are also called \emph{coisotropic
actions}. This is justified by the following discussion.

\begin{defini}
A submanifold $N$ of a symplectic manifold $(M,\omega)$ is said to
be \emph{coisotropic} if and only if for every $x\in N$, $(T_x
N)^{\perp_{\omega}} \subset T_x N$. In particular a $G$-action on
$M$ is called coisotropic if and only if there exists an open
dense subset $U\subset M$ with $Gx$ coisotropic for every $x\in
U$.
\end{defini}
\begin{lemma}
$M$ is a $G$ multiplicity-free space if and only the $G$-action on
$M$ is coisotropic.
\end{lemma}
\begin{proof}
First of all we note the following easy fact: if $f\in
C^{\infty}(M)^G$ then for every $\xi \in \mf g$ we have $\{f,
f_{\xi} \}=0,$ where $f_{\xi}$ is defined by
$f_{\xi}(x)=\mu(x)(\xi).$

Assume that a generic orbit $Gx$ is coisotropic.

Let $f,g \in C^{\infty}(M)^G.$ Since  $\{ f, f_{\xi} \}=\{g, f_\xi \}=0$ 
we have $X_f, X_g \in (T_x Gx)^{\perp_{\omega}}
\subset T_x Gx,$ since $Gx$ is coisotropic, for every $x\in U.$
Hence 
\[
\{f,g \}(x) =\omega(X_f,X_g)=0, \ \forall x \in U,
\]
which implies $\{f,g\}=0$, since $U$ is an open dense subset.

Vice-versa, let $x\in M$ be a regular point. By the slice-theorem
there are functions $f_1,\ldots,f_k \in C^{\infty}(M)^G$ with
$df_1 \wedge \ldots \wedge df_k \neq 0$ in some neighborhood $W$
of $Gx$ and
\[
Gx=\{ y \in W : f_1(y)=\ldots=f_k(y)=0 \}.
\]
From $\{f_i,f_j \}=0$ one may deduce that $X_{f_i} \in T_x Gx.$
Therefore, since $X_{f_i} \in (T_x Gx)^{\perp_{\omega}},$ $i=1,
\ldots, k$  and $X_{f_1},\ldots ,X_{f_k}$ are independent in $W$,
we have that $Gx$ is a coisotropic submanifold.
\end{proof}
\begin{rem}
Our proof is essentially one given in \cite{HW}. However in
\cite{HW} the authors assumed that $G$ is a compact Lie group,
their proof works also when the $G$-action is a proper action. 
\end{rem}
It is standard that given a $G$-orbit $Gx=G/G_x$, study
the slice representation, i.e. the linear representation 
of $G_x$ induced from the $G$-action 
on $T_x M / T_x Gx$. In \cite{HW} p. 274, as a consequence of the
arguments used in the Restriction Lemma, it was proved that given a
complex orbit $Gp=G/G_p$ then $G$ acts coisotropically on $M$ if and only
if $G_p$ acts coisotropically on the slice, 
whenever $M$ is a compact K\"ahler manifold and $G$ is a subgroup 
of its full isometric group. Here, we give the same result in symplectic
context.
\begin{prop}
Let $(M,\omega)$ be a symplectic manifold and let $G$ be a Lie
group which acts properly and in a Hamiltonian fashion on $M$ with
moment map $\mu: M \lra \mf g^*$. 
Let $Gx$ be a symplectic orbit. If $M$ is a $G$
multiplicity-free space then the slice representation is a
multiplicity-free representation. Moreover, if $G$ is compact the
vice-versa holds as well.
\end{prop}
\begin{proof}
It follows from symplectic slice, see \cite{lb}, \cite{or}
\cite{sl}.

There exists a neighborhood of the orbit $Gx$ which is
$G$-e\-qui\-va\-riantly symplectomorphic to a neighborhood of the zero
section of the symplectic manifold $ (Y=G \times_{G_x}( (\mf
g_{\beta} / \mf g_x )^* \oplus V), \tau)$, see \cite{lb},
\cite{sl} for an explicit description of the symplectic form
$\tau$,  with a $G$-moment map $\mu$ given by
\[
\mu([g,m,v])=Ad(g)(\beta+j(m)+ i(\mu_{V}(v))),
\]
where $\beta=\mu(x)$, $i: \mf g^*_x \hookrightarrow \mf g^*$ is
the transpose of the projection $p:\mf g \lra \mf g_x$, $j:( \mf
g_{\beta} / \mf g_x )^* \hookrightarrow \mf g^o_x$ ($\mf g^o_x$ is
the annihilator of $\mf g_x$ in $\mf g^*$) is defined by a choice
of a $G_x$-equivariant splitting and  finally $\mu_V$ is
the moment map of the $G_x$-action on the symplectic subspace $V$
of $(T_x Gx, \omega(x))$. Note that $V$ is isomorphic to the
quotient $((T_x Gx)^{\perp_{\omega}} / ((T_x Gx)^{\perp_{\omega}}
\cap T_x Gx))$. In the sequel we denote by $\omega_V=
\omega(x)_{|_{V}}$.

Since $Gx$ is symplectic, $ (Y=G \times_{G_x} V, \tau)$
and the moment map $\mu$ becomes
\[
\mu([g,m])=Ad^*(g)(\beta+ i(\mu_{V}(m))),
\]
Assume that $M$ is a $G$-multiplicity-free space and let
$y=[e,m]\in Y$ be such that $Gy$ is coisotropic.

Let $Y \in Ker {\rm d} (\mu_{V})_m$. Noting ${\rm
d}\mu_{[e,m]}(0,Y)={\rm d}( \mu_V)_m (Y)=0$, which implies, from
(\ref{nucleo}), $Y\in (T_y Gy)^{\perp_{\omega}} \subset T_y Gy$.
This claims $Y\in T_y Gy \cap V= T_m G_x m$, i.e the slice 
representation is multiplicity-free.

Assume that $G$ is a compact Lie group. It is well known
\[
\rm{cohom}(G,M)=\rm{cohom}(G_x ,V),
\]
which follows from the classical slice theorem for proper actions \cite{Path},
and  rk($G$)=rk($G_x$),
since $Gx$ is a symplectic manifold, where cohom($G$,$M$) denotes
the codimension of a principal orbit and for every compact group
$K$, rk($K$) denotes the rank, namely the dimension of the maximal
torus. If $G_x$ acts multiplicity-free on $V$ then
cohom($G_x,V)$=rk($G_x$)-rk($G_{{\rm princ}}),$ see \cite{HW},
where $G_{{\rm princ}}$ is the principal isotropy subgroup of the
action, which implies that 
\[
{\rm cohom}(G,M)=rk(G)-rk(G_{{\rm princ}}).
\]
Therefore, from Theorem $3$ \cite{HW} p.
$269$, we get $G$ acts multiplicity-free on $M$. 
\end{proof}
\begin{cor}\label{hnc}
Let $M$ be an irreducible Hermitian irreducible symmetric space of non compact
type. Let $G$ be a compact group which acts  in a Hamiltonian
fashion on $M$. Then $G$ acts multiplicity-free on $M$ if and only
if it acts multiplicity-free on $M^*$, the corresponding irreducible Hermitian
symmetric space of compact type.
\end{cor}
\begin{proof}
Since $G$ is compact it has a fixed point $x\in M$, from a Theorem
of Cartan, see \cite{he}. Hence $G$ acts multiplicity-free on $M$
if and only if $G$ acts multiplicity-free on the tangent space at
$x$ if and only if it acts multiplicity free on $M^*$.
\end{proof}
\begin{rem}
Corollary \ref{hnc} classifies completely the  compact Lie groups acting in
a Hamiltonian fashion and coisotropically on 
the irreducible Hermitian symmetric spaces of noncompact type,
due the results proved in \cite{BG}, \cite{Bi}, \cite{PTh}.
\end{rem}
\section{Equivalence Theorems for multiplicity-free action}
From now on we assume that $(M,\omega)$ is a connected symplectic
manifold acted on properly and in a Hamiltonian fashion by a connected
Lie group $G$. We denote by $\mu:M \lra \mf g^*$ the corresponding moment
map for the $G$-action on $M$.

Let $\alpha \in \gd$. We define the corresponding reduced space
\[
\sr{\alpha}=\mu^{-1}(G\alpha)/ G,
\]
to be the topological quotient of the subset $\mu^{-1}(G \alpha)$
of $M$ by the action of $G.$  It is well known, see \cite{acg},
\cite{lb}, \cite{sl}, that
the reduced space $\sr{\alpha}$ is a union of a symplectic
manifolds and it can be endowed  with a Poisson structure which arise from 
Poisson structure on the orbits space. Hence $\sr{\alpha}$ is a
symplectic stratified space and the manifolds which decompose
$\sr{\alpha}$ are called pieces.

Here we analyze the case when $G=G_1 \times G_2,$ where 
$G_i, \ i=1,2$ are closed connected subgroup of $G$. We assume also that
the coadjoint orbits of $G, G_1$ and $G_2$ are locally closed spaces.

Obviously  $\mf g^*= \mf g_1^* \oplus \mf g_2^*$ and the 
moment map $\mu=\mu_1 + \mu_2$, where $\mu_i$ is the corresponding 
moment map for the $G_i$-action on $M,$ $i=1,2. $ 
Since $\mu$ is $G$-equivariant, we
have  that $\mu_1$ is invariant under $G_2$ and $\mu_2$ is
invariant under $G_1$.

Let $\alpha=\alpha_1 + \alpha_2.$ The $G_1$-action  on the
pieces of  the reduced space $\sr{\alpha_2}=\mu_2^{-1}(G_2
\alpha_2)/G_2$ is symplectic. These moment maps on the pieces fit
together to form an application 
\[
\mu_{1,2}: \sr{\alpha_{2}} \lra
\mf g^*_1 ,
\]
such that $\mu_{1,2} =\mu_1 \circ \pi_2,$ where
$\pi_2$ is the projection on $\sr{\alpha_2}.$ Clearly  $G_2$ acts
on the reduced space 
$\sr{\alpha_1}=\mu_1^{-1}(G_1 \alpha_1) /G_1$ with moment map
\[
\mu_{2,1}: \sr{\alpha_{1}} \lra \mf g^*_2
\] 
such that $\mu_{2,1}=\mu_2 \circ \pi_1,$ where $\pi_1$ is the projection on
$\sr{\alpha_1}.$

We introduce the notion of multiplicity-free space for the reduced
space. Indeed, we say that the $G_1$-action on $\sr{\alpha_2}$ is
multiplicity-free if the ring of $G_1$-invariant functions
of  $\sr{\alpha_{2}}$, \cite{acg}, 
is a commutative Poisson algebra. 

We may also define the reduced space with respect the $G_1$-action 
on $\sr{\alpha_2}$ to be
\[
(M_{\alpha_2})_{\alpha_1}= \mu_{1,2}^{-1}(G_1 \alpha_1 )/G_1.
\]
The same definition holds for the $G_2$-action on $\sr{\alpha_1}.$

Now, we shall give a criterion for a $G$-action to be a
multiplicity-free action. We begin with the following lemma.
\begin{lemma}
Let $(M,\omega)$ be a symplectic manifold  and let $G=G_1 \times G_2$
be a connected Lie group which acts in a Hamiltonian fashion on
$M$. Let $\alpha_1 \in \mf g_1^*.$ Then  the $G_2$-action on
$\sr{\alpha_1}$ is multiplicity-free if and only if for every
$\alpha_2 \in \mf g_2^*$ the reduced space
$(\sr{\alpha_1})_{\alpha_2}$ are points
\end{lemma}
\begin{proof}
In the sequel we always refer to \cite{acg} and \cite{lb}. 

Let $\mu_2 : M  \lra \mf g_1^*$ be the moment map of the $G_1$-action
and let $\mu_{2,1}: \sr{\alpha_1} \lra \mf g_2^*$ be the
corresponding moment map of the $G_2$-action on the reduced space.
We recall that the smooth function on the reduced spaces are
defined by
\[
C^{\infty}(\sr{\alpha_1})= C^{\infty}(M)^{G_1} |_{\mu_1^{-1} (G_1
\alpha_1)}
\]
and the reduced space is a locally finite union of symplectic
manifolds (symplectic pieces) which are the following manifolds.

Let $H$ be a subgroup of $G.$ The set $M^{(H)}$ of orbit of type
$H,$ i.e. the set of points which orbits are isomorphic to $G/H,$ is a
submanifold of $M$ \cite{Path}. The set 
$(\mu_1^{-1}
(G_1 \alpha_1 ) \cap  M^{(H)})
$
is a submanifold of constant rank and the quotient 
\[
(
\sr{\alpha_1})^{(H)}=(\mu_1^{-1} (G_1 \alpha_1 ) \cap
M^{(H)})/G_1,
\]
is a symplectic manifold which inclusion $(
\sr{\alpha_1})^{(H)} \hookrightarrow \sr{\alpha_1}$ is a Poisson
map \cite{lb}.

The $G_2$-action preserves $(\sr{\alpha_1})^{(H)},$ and the following
topological space
\[
( ( \sr{\alpha_1})^{(H)} \cap \mu_{2,1}^{-1}(G_2 \alpha_2))
/G_2=\cup_{i \in I} S_i
\]
is a stratified symplectic manifold which restrictions map
\[
r^{H}_{\alpha_2}:C^{\infty}(( \sr{\alpha_1})^{(H)} \cap
\mu_{2,1}^{-1}(G_1 \alpha_1) )^{G_2} \lra C^{\infty}(S_i)
\]
are Poisson and surjective. Therefore, if
$C^{\infty}(\sr{\alpha_1})^{G_2}$ is abelian, the algebra
$C^{\infty}(S_i),\ i\in I$ must be abelian, and $S_i$ must be
discrete and therefore a points.

Vice-versa, if all reduced spaces are points then
$r^{H}_{\alpha_2}(\{f,g \})=0$ for all $\alpha_2 \in g_2^*,$ and
every $H$ subgroup of $G_1,$ so that $\{ f,g \}=0.$
\end{proof}
\begin{thm} \label{ghigi}
Let $(M, \omega)$ be a symplectic manifold with a proper and
Hamiltonian action of a connected Lie group $G=G_1 \times G_2.$
Assume also that $G_i,\ i=1,2$ are closed connected Lie group and
the coadjoint orbits of $G, G_1$ and  $G_2$ are locally closed. 
Hence $M$ is a $G$ multiplicity-free space if and only
if for every $\alpha=\alpha_1+ \alpha_2 \in \gd$ $\sr{\alpha_1}$
is a $G_1$ multiplicity-free space if and only if $\sr{\alpha_2}$
is a $G_2$ multiplicity-free space.
\end{thm}
\begin{proof}
It follows immediately from the above result. Indeed, it is easy
to check that $\sr{\alpha}=\mu^{-1}(G\alpha) /G$ is homemorphic to
$(M_{\alpha_1})_{\alpha_2}$ or equivalently is homemorphic to
$(M_{\alpha_2})_{\alpha_1};$ the homeomorphism is given by the
natural application
\[
( M_{\alpha_1} )_{\alpha_2} \longrightarrow \sr{\alpha}, \ \ \
[[x]] \longrightarrow [x]
\]
which preserves the symplectic pieces, concluding the prove.
\end{proof}
Theorem \ref{ghigi} is not difficult to prove. However, from Theorem 
\ref{ghigi}, we may deduce some interesting facts.
\begin{prop} \label{restriction}
Let $N$ be closed $G-$invariant symplectic submanifold of $M.$ If
$M$ is a $G$ multiplicity-free space then so is $N$.
\end{prop}
\begin{proof}
$G$ acts on $N$ in a Hamiltonian fashion with moment map
$\overline{\mu}: N \lra \mf g^* ,$ $\overline{\mu}(x) =\mu(x),$
which is the restriction of $\mu$ on $N$.  In particular, for
every $\alpha \in \mf g^*$ the reduced space $N_{\alpha} \subset
\sr{\alpha},$ which implies that the topological space
$N_{\alpha}$ are points.
\end{proof}
\begin{cor}
If $G$ is a compact Lie group acting multiplicity-free on $M$ then
\[
M^G:=\{x \in M: Gx=x \},
\]
must be a finite set.
\end{cor}
Another interesting application of Theorem \ref{ghigi} is the
following result.
\begin{prop} \label{po}
Let $(M,\omega)$ be a symplectic manifold with a Hamiltonian circle
action. Let $K$ be a connected Lie group which acts properly and
in a  Hamiltonian fashion on $M.$ Assume also that the $K$-action
commutes with the circle action.
If $M$ is a $K$ multiplicity-free space then so is the $K$-action induced
on any symplectic cut.
\end{prop}
\begin{proof}
We briefly describe the symplectic cut, see \cite{lerman} and \cite{gbl} 
for more details.

Let $(M,\omega)$ be a symplectic manifold with a Hamiltonian action
of a one-dimensional torus $\To^1$ 
with moment map $\phi: M \lra \R.$ We consider the
symplectic manifold $N=M \times \C$, equipped with the symplectic
form $\omega - \frac{i}{2} dz \wedge d \overline{z}.$ 
$\To^1$ acts on $N$ with its product action and this action 
is a Hamiltonian action 
with moment map $\psi(p,z)=\phi(p)+ |z|^2.$  The reduced space
$\mk= \psi^{-1}(\lambda)/ S^1,$ $\lambda \in \R$ is the symplectic
cut.

The $K-$action on $M \times \C$ is given by $k(m,z)=(km,z).$ Since
the $K-$action commutes with the $\To^1$-action, it 
induces a Hamiltonian action on the symplectic cut  with
moment map $\overline{\mu}([x,z])=\mu(x)$ where $\mu$ is the
moment map of the $K-$action on $M.$ 
Note that $K \times \To^1$-action is multiplicity-free on $M\times \C$ if
and only if the $K$-action is on $M$. Therefore, if $K$ acts multiplicity-free
on $M$,
from Theorem \ref{ghigi}, $K$ acts multiplicity-free on the symplectic cut. 
\end{proof}

Let $H$ be a compact subgroup of $G$ and let $N(H)$ be its
normalizer in $G$. It is well-known that the Lie group $L=N(H)/H$  
acts freely and properly
on the submanifold $M^{H}:=\{ x\in M: G_x=H \}$ \cite{Path}. Moreover, since
$T_x M^H = (T_x M  )^H$, $M^H$ is  symplectic.

In \cite{lb} it was proved that
$L$ acts in a Hamiltonian fashion on $M^{H}$, the dual of the Lie algebra of
$L$ is naturally isomorphic to the subspace $(\mf h^o)^H$ of the
$H$-fixed vectors in the annihilator of $\mf h$=Lie($H$) in 
$\mf g^*$. Furthermore, given $\alpha=\mu(x)$, where $x\in M^{H}$, then
\[
G\mu^{-1}(\alpha) \cap M^{H}/ G \cong (M^{H})_{\alpha_o },
\]
$\cong$ means symplectically diffeomorphic,
where $\alpha_o=\pi(\alpha)$ and  
$\pi: (\mf h^o)^{H} \lra l^*$ be the natural projection. This proves that
if $M$ is a $G$ multiplicity-free space then, from Theorem \ref{ghigi},
$M^H$ is a $L$ multiplicity-free space. Hence, we have the following
result.
\begin{prop} \label{oba}
Let $H$ be a compact subgroup of $G$. If $G$ acts coisotropically
on $M$ then $L$ acts coisotropically on $M^{H}$.
\end{prop}
Next, we claim the reduction principle for a multiplicity-free action.
\begin{prop}(\textbf{Reduction principle})
Let $G$ be a connected Lie group acting properly on a connected
symplectic manifold $M$. Let $H$ be a principal isotropy for the
$G$-action. Then $G$ acts coisotropically on $M$ if and only if
$L=N(H)/H$ acts coisotropically on $M^{H}$.
\end{prop}
\begin{proof}
Since the $G$-action is proper and preserves $\omega$, 
there exists a $G$-invariant almost complex structure $J$, i.e. $J:TM \lra TM$
be such that $J^2 =-Id$, adapted to $\omega$, i.e. 
$\omega(J\cdot, J\cdot)=\omega(\cdot, \cdot )$ and $g=\omega(\cdot, J \cdot)$
is a Riemannian metric, see
\cite{lb}. 

Now let $H$ be a principal isotropy and let $L=N(H)/H$. It is
well-known that 
\[
M^{H} \cong N(H)/H \times (T_{x}Gx)^{\perp_g },
\] 
see \cite{Path}, which implies $T_x Lx= (T_x Gx )^H$. 

Since
$(T_x Gx )^{\perp_{\omega}}=J(( T_x Gx )^{\perp_g})$ and 
$(T_x Gx)^{\perp_g}\subset (T_x M )^H$, recall that $H$ acts trivially on the 
slice due the fact that $Gx$ is a principal orbit, we get that 
\[
J( (T_x Gx )^{\perp_g }) \subset T_x Gx \Leftrightarrow 
J( (T_x Gx )^{\perp_g }) \subset (T_x Gx )^H .
\] 
Therefore, recall that $(T_x Lx )^{{\perp_g }_{|_{T_x M^H }}}=(T_x Gx )^{\perp_g }$
since $Gx$ is principal, we have that $Gx$ is coisotropic in $M$ if and only if
$Lx$ is in $M^H$.
\end{proof}
We conclude this section 
giving the Equivariant mapping lemma, see \cite{HW}, in a
symplectic context.
\begin{prop} \label{equi}
Let $(M,\omega)$ and $(N, \omega_o)$ be connected symplectic
manifolds and $G$ be a connected Lie group acting on both
manifolds properly, and in a Hamiltonian  fashion. Let $\phi: M
\lra N$ be a smooth surjective $G-$equivariant map with
${\rm rank} \phi=\dim N.$ Assume that for every $p \in M,$  $Ker d
\phi_p$ is a symplectic subspace and 
\[
d\phi_p : ((Ker
\phi)^{\perp_{\omega}}, \omega)  \lra (T_p N, \omega_o)
\] 
is a symplectomorphism. If $M$ is a $G$ multiplicity-free
space then so
is $N$.
\end{prop}
\begin{proof}
Let $f \in C^{\infty}(N)^G .$ The function $\tilde{f}= f \circ
\phi$ is a $G$-invariant function of $M.$ Take the vector field
$X_f$ such that $df= i_{X_f } \omega_o.$ By assumption the vector
field $\tilde X \in (Ker \phi)^{\perp_{\o}}$ such that $d\phi(
\tilde X )= X_f$ is the  symplectic gradient of $\tilde f.$ Hence,
given $f,g \in C^{\infty}(N)^G$ there exist $\tilde f,$ $\tilde g
\in C^{\infty}(M)^G$ such that $\{ f, g \}(\phi(x) )= \{\tilde f ,
\tilde g \}(x)$ which conclude our proof.
\end{proof}
%
%
\section{multiplicity-free spaces in K\"ahler geometry} 
\label{pussavia}
In the sequel we shall assume that $M$ is a compact K\"ahler manifold
and $G$ is a closed subgroup of its full isometry group 
acting in a Hamiltonian fashion
on $M$. Note that this action is automatically holomorphic by  
a Theorem of Konstant (see \cite{KN} p.242). 

In \cite{Pu}
it was introduced the {\em homogeneity rank} of $(G,M)$ as the
following integer
$$
{\rm homrk}(G,M)= {\rm rk}(G)- {\rm rk} (G_{{\rm princ}})- {\rm
cohom}(G,M),
$$
where $G_{{\rm princ}}$ is the principal isotropy subgroup of the
action, the integer cohom$(G,M)$ is the codimension of the
principal orbit and, for a compact Lie group $H,$  {\rm rk}$(H)$
denotes the rank, namely the dimension of the maximal torus.

In \cite{HW} it was proved that a group $G$ acts multiplicity-free
on $M$ if and only if the homogeneity rank vanishes. 

Our purpose is to make a 
connection between homogeneity rank of $(G,M)$ and homogeneity rank
of $(G, \mr)$, where
$\mr$ is the reduced space obtained from a torus action. 

Let $K$ be a semisimple compact Lie subgroup of $G$ and let $\To^k$ be a
$k$-dimensional connected torus which centralizes $K$ in $G$, i.e. $\To^k
\subset C_G (K)$. In the sequel we denote by
\[
\phi:M \lra \mf k \oplus  \mf t_k,
\]
where $\mf t_k=$Lie$(\To^k)$, the moment map of the $\To^k \cdot K$-action
on $M$ and with $\mu$, respectively with $\psi$, 
the moment map of the $K$-action on $M$, respectively 
a moment map of the $\To^k$-action on $M$.

Let $\lambda \in \mf t_k$ be such that
$\To^k$ acts freely on $\psi^{-1}(\lambda)$. The symplectic reduction
\[
(M_{\lambda}= \psi^{-1}(\lambda) / \To^k, \omega_{\lambda}),
\] 
is a symplectic manifold and $\omega_{\lambda}$ satisfies
\[
\pi^*(\omega_{\lambda})= i^*(\omega),
\]
where $\pi$ is the natural projection $\psi^{-1}(\lambda) 
\stackrel{\pi}{\lra} \mr$ and
$i$ is the inclusion $\mk \hookrightarrow M$,
\cite{Ca}, \cite{smd}.
Since $K$ commutes with
$\To^k$, $K$ acts on $M_{\lambda}$ in a Hamiltonian fashion with moment map
\[
\overline{\mu}: M_{\lambda} \lra \mf k , \ \overline\mu ([x])= \mu(x).
\]
Indeed, It is easy to see that $\overline{\mu}$ is $K-$equivariant.
Hence the problem is then restricted
to verify that for every $Z \in \mf k$  we have
$d \overline{\psi}^{Z} = i_{ \tilde {Z}^\# } \omr,$
where $\tilde{Z}^\#$ is the vector field on $\mr$ generated by the
one parameter subgroup $\exp(tZ).$

Let $X \in T_{[x]} \mr$ and
let $\tilde{X} \in T_{x} \mu^{-1}(\lambda)$ such that
$\pi_* (\tilde{X}) = X.$  Since $\pi_* (Z^\#)=\tilde{Z}^\#$, where $Z^\#$
is the Killing field induced from $Z$ in $M$,
it follows
$$
d \overline{\psi}^Z (X) = d \psi^Z (\tilde{X})= 
i_{Z^\#} \omega (\tilde{X})= \pi^* \omr (Z^\#, \tilde{X})
=  i_{\tilde{Z}^\#} \omr (X),
$$
thus $K$ acts in a Hamiltonian fashion on $\mk.$

Let $[p]\in M_{\lambda}$. It is easy to see that $k[p]=[p]$
if and only if there exists $r(k)\in \To^k$ such that $kp=r(k)p$,
which is unique  since $\To^k$ acts freely on
$\psi^{-1}(\lambda)$. This means that the following application
\begin{equation} \label{pq}
K_{[p]} \stackrel{F}{\lra} (\To^k \cdot K)_p,\ F(k)=k r(k)^{-1},
\end{equation} 
is a covering map, due the fact that $K$ is semisimple. Hence
\begin{equation} \label{pw}
\dim K[p]= \dim (\To^k \cdot K)p - \dim \To^k.
\end{equation}
Since $M$ is a compact manifold, we
may extend the $\To^k$-action to a holomorphic action of
$(\C^*)^n$ which commutes with the $K-$action. 
This can be easily deduced from the following easy   
fact: let $X,Y$ be holomorphic fields. If $[X,Y]=0$ then
$[X,J(Y)]=0$, since 
$[X,J(Y)]=J([X,Y])=0$, due the fact that $M$ is K\"ahler.  
In particular the infinitesimal generatores of the $K-$action commute with ones
of the $(\C^* )^n -$action, proving that the two action commute as well.

The set $(\C^* )^n \cdot
\psi^{-1}(\lambda)$ is an open subset. Indeed, for every $p \in
\psi^{-1}(\lambda)$, denoting with $\mf z_p$ the vector subspace of
$T_p M$ spanned by the infinitesimal generator of the
$\To^k$-action on $M$, we have $T_p \psi^{-1}(\lambda) \oplus
J(\mf z_p)=T_p M$, since $\lambda$ is a regular value, 
which implies our affirmation. In particular the
open subset $(\C^*)^n \cdot \psi^{-1}(\lambda)$  contains regular
points. Hence there exists an element
\[
q=\rho_1 \cdots \rho_n \exp(i \theta_1) \cdots \exp(i \theta_n)p
=\rho \exp(i \Theta)p \in (\C^*)^n \cdot \psi^{-1}(\lambda),
\]
such that the orbit $( \To^k\cdot K) q$ is a principal orbit. 
Since $K$ commutes with $(\C^* )^n$,
we get that $(\To^k \cdot K)_p = (\To^k \cdot K )_q$ which means
that $p$, which lies in $\psi^{-1}(\lambda)$, is a regular point. 
Therefore, from (\ref{pq}) we deduce that $K[p]$ is a principal orbit and
from (\ref {pw}), we get 
\begin{equation}
{\rm homrk} (K, \mr)= {\rm homrk}(\To^k \cdot K, M),
\end{equation}
which proves the following result.
\begin{prop} \label{red}
Let $G$ be a connected compact Lie group
acting isometrically and in a Hamiltonian 
fashion on a compact K\"ahler manifold $M$. Let $K$ be a compact semisimple Lie
group of $G$ which centralizer in $G$ contains a $k$-dimensional connected
torus $\To^k$. Let 
$\lambda \in \mf t_k$ be such  that 
$\To^k$ acts freely on $\psi^{-1}(\lambda)$, where $\psi$ is a moment map of 
$\To^k$-action on $M$. Then the $(\To^k \cdot K)$-action is 
coisotropic on $M$ if and only if the $K$-action is on  
$M_{\lambda}=\psi^{-1}(\lambda)/ \To^k$.
\end{prop}
If we consider a one-dimensional torus $\To^1$ we may 
investigate the $K$-action on the K\"ahler cut $\mk$ obtained from the 
$\To^1$-action. Here we only assume that the $K$-action commutes with the 
$\To^1$-action. It is easy to check that $K_{[v,z]}=K_v$ when 
$z \neq 0$. Since $\{[(v,z]\in \mk : \ z\neq 0 \}$ is an open subset,
one may deduce that
\begin{equation} \label{k2} {\rm homrk}(K,M)={\rm
homrk}(K,\mk).
\end{equation}
Hence $K$ acts coisotropically on $M$ if and only if $K$ acts on $\mk$
which proves the following result
\begin{prop} \label{cut}
Let $G$ be a connected compact Lie group
acting isometrically and in a Hamiltonian 
fashion on a compact K\"ahler manifold $M$. Let $K$ be a compact Lie
group of $G$ whose centralizer in $G$ contains a one-dimensional
torus $\To^1$. Let 
$\lambda \in \mf t_1$ be such  that 
$\To^1$ acts freely on $\psi^{-1}(\lambda)$, where $\psi$ is a moment map of 
$\To^1$-action on $M$. Then $K$-action is 
coisotropic on $M$ if and only if the $K$-action on  
the K\"ahler cut $\mr$ is.
\end{prop}
\begin{exe}
Let $\omega=  \sqrt{-1} \sum_{i=1}^{n+1} d z_i \wedge d
\overline{z}_i$ be a K\"ahler form on $\C^{n+1}.$ Consider the
following $S^1 -$action on $(\C^{n+1}, \omega):$
$$
t \in S^1 \mapsto \theta_t={\rm multiplication\ by\ }e^{it}.
$$
The action is Hamiltonian with moment map $\mu(z)= | z |^2+$
constant. If we choose the constant to be $-1,$ then
$\mu^{-1}(0)=S^{2n+1}$ is the unit sphere on which $S^1$ acts
freely and the K\"ahler reduction $\mu^{-1}(0)/S^1$ is just $\P^n
( \C)=\SU(n+1) / \SUr{1}{n}.$ Therefore, by Proposition \ref{red}, 
a compact connected Lie subgroup $K$ 
of $\SU (n+1)$ acts multiplicity-free on $\P^n(\C)$ 
if and only if $S^1 \cdot K$ acts multiplicity-free on
$\C^n.$ Kac \cite{Kac} and Benson and Ratclif \cite{BR} have given
the complete classification of linear coisotropic actions, from which one
has the full classification of coisotropic actions on $\P^n (\C).$

If we consider the cut of $\C^{n+1}$ at $\lambda>0$, with respect the above
${\rm S}^1$-action, we obtain, see \cite{gbl}, $\P^{n+1}(\C)$, with $\lambda$
times the Fubini-Study metric. Hence $G\subset \SU(n+1)$ acts
coisotropically on $\P^{n+1}(\C)$ if and only if it acts coisotropically on
$\C^{n+1}$.
\end{exe}
\section{Multiplicity-free actions on compact non Hermitian 
symmetric  spaces}  \label{quo}
Let $\To^1$ acting on $\P^n (\C),$ as
$$
(t, [z_o, \ldots, z_n])  \longrightarrow [z_o,\cdots ,tz_n] .
$$
This is a Hamiltonian action with moment map
$$
\phi ([z_o, \ldots, z_n] )= \frac{1}{2} \frac{\pr z_n \pr^2}{\pr
z_0 \pr^2 + \cdots + \pr z_n \pr^2}  \ .
$$
Note that $\phi([0, \ldots , 1])$ is the maximum value
of $\phi$ and $\phi^{-1} (\frac{1}{2})=[0, \ldots ,
1].$ Hence (see \cite{gbl} page 5) if $\lambda=\frac{1}{2}-
\epsilon,$ $\epsilon \cong 0,$ then the K\"ahler cut 
$\P^n (\C)^{\lambda}$ is the blow up of $\P^n (\C)$ at
$[0, \ldots , 1].$ 

Let $\To^n$ be a torus acting on $\P^n (\C)$ as follows
\[
(t_1,\ldots,t_n)([z_o, \ldots,z_n])=
(t_1 z_o, t_2 z_1, \ldots, t_n z_{n-1}, z_n]) .
\]
This action is Hamiltonian and the principal orbits are Lagrangian; therefore
$\To^n$ acts co\-is\-o\-tro\-pi\-cal\-ly 
on $\P^n (\C)$. Since $\To^n $-action commutes with the
above $\To^1$-action, from Proposition \ref{cut}, we get $\To^n$ acts
coisotropically on the blow-up at one point of $\P^n (\C)$.

We may generalize the above procedure as follows.

Let $G$ be a connected compact Lie group acting 
coistropically on a compact K\"ahler manifold. It is well-known that, 
see \cite{he}, $\mf g= \mf{z(g)}\oplus [\mf g, \mf g]$, and if we denote by
$G_{{\rm ss}}$ the semisimple connected compact Lie group whose Lie algebra is 
$[\mf g , \mf g]$, then
\[
G= Z(G) \cdot G_{{\rm ss}}.
\]
From Proposition \ref{red}, if $G_{ss}$ acts coisotropically on $M$, then so is
the $G_{ss}$-action induced on $\mk$, the reduced space obtained from
$\To^k \subset Z(G)$.

Let $\To^1$ be a one-dimensional torus which lies on $Z(G)$. If $K\subset G$ 
is a compact Lie group acting coisotropically on $M$ then from
Proposition \ref{cut} $K$ acts coisotropically on the K\"ahler cut, obtained
from the $\To^1$-action on $M$. In particular, see \cite{gbl}, if
$\lambda_o$ is a maximum for the moment map of the $\To^1$-action then
$M^{\lambda}$, $\lambda=\lambda_o - \epsilon$, $\epsilon \cong 0$, is the 
blow-up of $M$ along the complex submanifold 
$\psi^{-1}( \lambda_o)$, where $\psi$ is the corresponding 
moment map for the $\To^1$-action on $M$.

In \cite{BG},\cite{Bi}, \cite{PTh}, the complete classification
of coisotropic actions on irreducible Hermitian symmetric spaces of
compact type is given. Therefore, it is easy to construct examples using the 
above strategy. For example, the $\SU(n)$-action  on $\SO(2n)/\U(n)$
induces a coisotropic action on K\"ahler cut given by 
$Z(\U(n))$. More generally, if $M=L/P$ is an irreducible 
Hermitian symmetric space of compact type, 
then $Z(P)$ is a one-dimensional torus. Since the $P$-action on $M$ is 
coisotropic, see \cite{BG}, \cite{Bi}, \cite{PTh}, $P$ acts 
coisotropically
on the K\"ahler cut with respect  the $Z(P)$-action on $M$.   

\end{document}